# EFFECTIVE RESISTANCE OF RANDOM TREES

By Louigi Addario-Berry,[1] Nicolas Broutin and Gábor Lugosi[2]

*Université de Montréal, INRIA Rocquencourt and
ICREA Pompeu Fabra University*

We investigate the effective resistance $R_n$ and conductance $C_n$ between the root and leaves of a binary tree of height $n$. In this electrical network, the resistance of each edge $e$ at distance $d$ from the root is defined by $r_e = 2^d X_e$ where the $X_e$ are i.i.d. positive random variables bounded away from zero and infinity. It is shown that $\mathbf{E} R_n = n \mathbf{E} X_e - (\mathbf{Var}(X_e)/\mathbf{E} X_e) \ln n + O(1)$ and $\mathbf{Var}(R_n) = O(1)$. Moreover, we establish sub-Gaussian tail bounds for $R_n$. We also discuss some possible extensions to supercritical Galton–Watson trees.

**1. Introduction.** Electric networks have been known to be closely related to random walks and their investigation often offers an elegant and effective way of studying properties of random walks. See Doyle and Snell [9] and Lyons and Peres [15] for very nice introductions to the subject. For the better understanding of certain random walks in random environments, it is natural to study random electric networks, that is, electric networks in which edges are equipped with independent random resistances. This model was studied by Benjamini and Rossignol [5] who considered the case of the cubic lattice $\mathbb{Z}^d$ where the resistance of each edge is an independent copy of a Bernoulli random variable. Using an inequality of Falik and Samorodnitsky [11], they proved that point-to-point effective resistance has submean variance in $\mathbb{Z}^2$, whereas the variance is of the order of the mean when $d \geq 3$. In this paper, we study the corresponding problem for binary trees.

An electric network is a locally finite connected graph $G = (V, E)$ with vertex set $V$ and edge set $E$ such that each edge $e \in E$ is equipped with a number $r_e \geq 0$ called *resistance*. (In this paper we only consider finite

Received January 2008; revised October 2008.
[1]Supported in part by an NSERC discovery grant.
[2]Supported by the Spanish Ministry of Science and Technology Grant MTM2006-05650 and PASCAL 2 Network of Excellence.

graphs.) Alternatively, an edge is associated with a conductance $c_e = 1/r_e$. The *effective resistance* between two disjoint sets of vertices $A, B \subset V$ is defined as follows: assign "voltage" (or potential) $U(u) = 1$ to each vertex $u \in A$ and $U(v) = 0$ for all $v \in B$. If $G$ is finite then the function $U$ can be extended, in a unique way, to all vertices in $V$ according to two basic laws given by *Ohm's law* and *Kirchhoff's node law*. In order to describe these laws, we need the notion of *current*. Given two vertices $u, v \in V$ joined by and edge $e \in E$, the current flowing from $u$ to $v$ is a real number $i(u, v)$. Ohm's law states that for each edge of the graph, $i(u,v)r_e = U(u) - U(v)$. Kirchhoff's node law postulates that for any vertex $u \notin A \cup B$, $\sum_{v \,:\, v \sim u} i(u, v) = 0$. (For the proof that these two laws uniquely determine the function $U : V \to [0, 1]$, see [9] or [15].) Now the *effective conductance* between the vertex sets $A$ and $B$ is defined as the total current flowing into the network, that is,

$$C(A \leftrightarrow B) = \sum_{u \in A} \sum_{v \,:\, v \sim u} i(u, v).$$

The effective resistance between $A$ and $B$ is $R(A \leftrightarrow B) = 1/C(A \leftrightarrow B)$.

Several useful tricks of network reduction are known that help simplify resistance calculations. Since in this paper we focus on trees, it suffices to recall two of the simplest rules. One of them states that two resistors in series are equivalent to a single resistor whose resistance is the sum of the original resistances. The other rule states that two conductors in parallel are equivalent to a single conductor whose conductance is the sum of the original conductances. Apart from these two simple rules, a formula called *Thomson's principle* will be useful for our purposes. Thomson's principle gives an explicit expression for the effective resistance. It states that

$$(1) \qquad R(A \leftrightarrow B) = \inf_{\Theta \in F} \sum_{e \in E} r_e \Theta(e)^2,$$

where the infimum is taken over the set $F$ of all *unit flows*. A unit flow is a function $\Theta$ over the set of oriented edges $\{(u, v) : u \sim v\}$ which is antisymmetric [i.e., $\Theta(u, v) = -\Theta(v, u)$], satisfies $\sum_{v \,:\, v \sim u} \Theta(u, v) = 0$ for any vertex $u \notin A \cup B$, and has

$$\sum_{u \in A} \sum_{v \notin A \,:\, v \sim u} \Theta(u, v) = \sum_{v \in B} \sum_{u \notin B \,:\, u \sim v} \Theta(u, v) = 1.$$

It can be shown that the unique unit flow $\Theta^*$ which attains the above infimum is proportional to the current $i(u, v)$ (see, e.g., Doyle and Snell [9], page 50).

In this paper we consider the case of a complete infinite binary tree $T$ with root $r$. (All results carry over trivially to infinite $b$-ary trees for integers $b > 2$.) The *depth* $d(v)$ of a node $v$ in $T$ is the number of edges on the path from the root to $v$. We say that an edge $e$ has depth $d$ if there are $d$ edges on the path starting with edge $e$ and ending at the root. The resistance of an edge $e$ at depth $d$ is defined by $2^{d-1}X_e$ where the $X_e$ are independent copies of some strictly positive random variable with finite mean. This exponential weighting corresponds to the "critical" (with respect to transience/recurrence) case of the biased random walk in random environment obtained by traversing an edge $e$, starting from either endpoint, with probability proportional to its conductance (the inverse of its resistance). This type of exponential scaling of resistances was considered, for example, by Lyons [13]. He showed that in an infinite rooted tree with branching number $b$, if the (deterministic) resistance of an edge equals $\lambda^d$ then the effective resistance between the root and "infinity" is infinite if $\lambda > b$ and finite if $\lambda < b$. Thus, our choice of scaling corresponds to the critical case. Similar biased random walks have been studied in depth by Pemantle [17], Lyons [13] and Lyons and Pemantle [14], who beautifully characterize the type of such random walks in many situations. (However, our model does not quite fit within their framework, as the transition probabilities fail to satisfy a certain independence requirement.) Unfortunately, we do not immediately see how to translate our results into results about biased random walks or random walks in random environments. Also, it is likely that such results would only be new if we could also extend our results to the more general setting of Galton–Watson trees. For more background on the connection between effective resistance of networks and random walks, see Doyle and Snell [9], Lyons and Peres [15], Peres [18] or Soardi [19].

For a random network such as that described in the previous paragraph, interesting and nontrivial behavior emerges. Let $R_n$ be the effective resistance between the root $r$ and the set of vertices at depth $n$, and let $\mu$ and $\sigma^2$ be the mean and variance of $X_e$, respectively. The primary results of our paper are that as long as $X_e$ is bounded away from both zero and infinity,

$$\mathbf{E}R_n = \mu n - \frac{\sigma^2}{\mu}\ln n + O(1) \quad \text{and} \quad \mathbf{E}[|R_n - \mathbf{E}R_n|^q] = O(1) \qquad \text{for all } q \geq 1.$$

(These results are contained in Theorems 5 and 7.) We also derive correspondingly precise results about the conductance $C_n = 1/R_n$. Interestingly, in order to estimate the expected resistance, our main tool is a sharp upper bound for the variance of the conductance (and thereby for the variance of the resistance). Intuitively, concentration of the conductance implies that the behavior of the electric network is not very different from the one with deterministic resistances $2^d\mu$. Thus, Section 2 is devoted to the variance of the conductance $C_n$. In particular, we show that $\mathbf{Var}[C_n] = O(n^{-4})$. In Section 3 we derive the bounds for the expected resistance and conductance mentioned above. In Section 4 we establish sub-Gaussian tail bounds for $R_n$. The proof is based on Thomson's formula (1) and relies on an exponential concentration inequality due to Boucheron et al. [7].

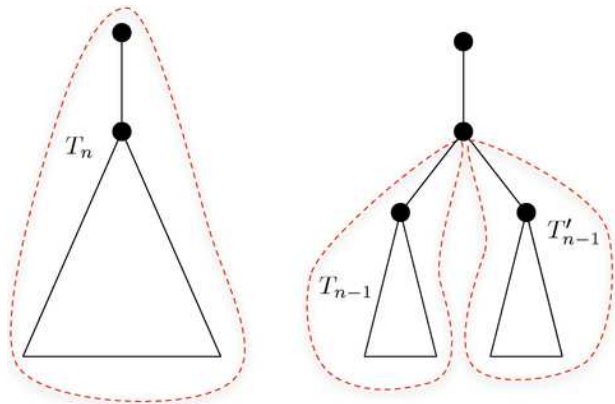

Fig. 1. *Rooting the binary tree rooted at an edge instead of a vertex simplifies the recursive decomposition.*

Finally, in Section 5 we briefly discuss the ways in which one might attempt to extend our results to supercritical Galton–Watson processes (in this case the appropriate scaling for the resistances is $[\mathbf{E}Z_1]^d$ for edges at depth $d$, where $Z_1$ is the number of offspring of the root). In the Galton–Watson setting, it makes sense to first condition on the tree, then study the conditional behavior of the effective resistance and of $\Theta^*$. In Section 5 we shall also observe that if the random variable $X_e$ is constant then the "scaled analogue" of Question 4.1 from Lyons, Pemantle and Peres [16] is easily answered; motivated by this, we suggest a more general question.

From this point on, we assume $X_e$ is any random variable taking values in some interval $[a,b]$ with $0 < a < b$. Most our arguments rely on a recursive decomposition of the tree. This decomposition is made easier by rooting the tree at an edge instead of at a vertex; this trick was also used in [16] to facilitate conductance computations. More precisely, for $n \geq 1$ we define the tree $T_n$ as follows: the root has one single child whose subtree is a complete binary tree with $n-1$ levels (so $2^{n-1}$ leaves). Then, $T_n$ decomposes exactly into a single edge connected (in series) to two independent copies of $T_{n-1}$ (in parallel) as shown in Figure 1. We let $R_n$ be the effective resistance of $T_n$ taken between the root and the leaves. Let $C_n = 1/R_n$ be the corresponding effective conductance, so in particular $R_1$ is distributed as $X$ and $C_1$ is distributed as $1/X$. The difference between $R_n$ and the effective resistance of the complete binary tree of height $n-1$ is at most $b$, so bounds on the moments of the former immediately imply corresponding bounds for the latter.

We close this introduction by noting that the results of Benjamini and Rossignol [5] are proved by adapting an argument first used by Benjamini, Kalai and Schramm [6] to prove submean variance bounds for first-passage percolation on $\mathbb{Z}^2$. Addario-Berry and Reed [1] have studied first-passage percolation on supercritical Galton–Watson processes; though their approach is entirely different from ours, their result is strikingly similar: under suitable assumptions on the edge lengths (which in their case are i.i.d.), the height of the first-passage percolation cluster of (weighted) diameter $n$ has expected value $\alpha n - \beta \ln n + O(1)$, for computable constants $\alpha$ and $\beta$, and has bounded variance. We are not sure whether this similarity is more than a coincidence.

**2. The variance of the conductance.** The purpose of this section is to derive an upper bound for the variance of the conductance $C_n$. We start by noticing that $R_n$ and $C_n$ admit the following scalings.

LEMMA 1. *When $a \leq X \leq b$, we have $an \leq R_n \leq bn$ and $1/b \leq nC_n \leq 1/a$, for all $n \geq 1$.*

The lemma follows from Rayleigh's monotonicity law (see [9], page 53) by bounding the resistance of $T_n$ between that of two deterministic networks in which the random variables either always take their minimum value $a$ or their maximum value $b$. We first derive a bound on $\mathbf{Var}[C_n]$. Using Chebyshev's inequality, this bound yields a quadratically decaying tail bound for $R_n$.

THEOREM 2. *There exists a constant $K$ depending only on $a$ and $b$ such that $\mathbf{Var}[C_1] \leq K, \mathbf{Var}[R_1] \leq K$ and for all integers $n \geq 2$,*

$$\mathbf{Var}[C_n] \leq K \cdot \left( \sum_{i=1}^{n-1} \frac{2^{1-i}}{(n-i)^4} + \frac{1}{2^{n-1}} \right) \leq \frac{2^{10} K}{n^4} \quad \text{and} \quad \mathbf{Var}[R_n] \leq K.$$

Our main tool in proving Theorem 2 is the Efron–Stein inequality, which provides an upper bound on the variance of functions of independent random variables.

THEOREM 3 (Efron and Stein [10], Steele [20]). *Let $Y_i$, $i \geq 1$, be independent random variables, and let $f : \mathbb{R}^n \mapsto \mathbb{R}$ be a measurable function of $n$ variables. Then,*

$$\mathbf{Var}[f(Y_1, \ldots, Y_n)]$$
$$\leq \frac{1}{2} \cdot \sum_{i=1}^{n} \mathbf{E}[(f(Y_1, \ldots, Y_i, \ldots, Y_n) - f(Y_1, \ldots, Y_i', \ldots, Y_n))^2],$$

*where $Y_i'$, $i \geq 0$, are independent copies of $Y_i$, $i \geq 0$.*

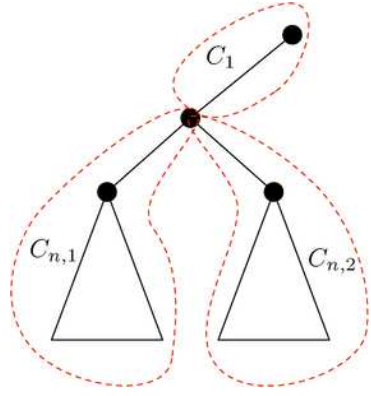

FIG. 2. *The decomposition of $T_n$ into three conductors $C_1$, $C_{n,1}$ and $C_{n,2}$ at the origin of our recurrence relations.*

PROOF OF THEOREM 2. It clearly suffices to treat the case $n \geq 2$. We decompose $T_n$ into three independent conductors $C_1$, $C_{n,1}$ and $C_{n,2}$ as shown in Figure 2. Then, $C_n$ is a function of these three independent random variables:

$$(2) \qquad C_n = \frac{C_1 \cdot (C_{n,1} + C_{n,2})}{C_1 + C_{n,1} + C_{n,2}}.$$

By the Efron–Stein inequality and the symmetry of $C_{n,1}$ and $C_{n,2}$, we have

$$(3) \quad \begin{aligned} \mathbf{Var}[C_n] &\leq \mathbf{E}\bigg[\bigg(\frac{C_1 \cdot (C_{n,1} + C_{n,2})}{C_1 + C_{n,1} + C_{n,2}} - \frac{C_1 \cdot (C'_{n,1} + C_{n,2})}{C_1 + C'_{n,1} + C_{n,2}}\bigg)^2\bigg] \\ &\quad + \frac{1}{2} \cdot \mathbf{E}\bigg[\bigg(\frac{C_1 \cdot (C_{n,1} + C_{n,2})}{C_1 + C_{n,1} + C_{n,2}} - \frac{C'_1 \cdot (C_{n,1} + C_{n,2})}{C'_1 + C_{n,1} + C_{n,2}}\bigg)^2\bigg], \end{aligned}$$

where $C'_{n,1}$ is an independent copy of $C_{n,1}$. Letting $\alpha = C_{n,1} + C_{n,2}$, the second term of (3) reduces to

$$\mathbf{E}\bigg[\bigg(\frac{C_1 \alpha}{C_1 + \alpha} - \frac{C'_1 \alpha}{C'_1 + \alpha}\bigg)^2\bigg].$$

Observe that, since $1/b \leq C_1, C'_1 \leq 1/a$, we have

$$\bigg|\frac{C_1 \alpha}{C_1 + \alpha} - \frac{C'_1 \alpha}{C'_1 + \alpha}\bigg| = \bigg|\frac{\alpha^2 (C_1 - C'_1)}{(C_1 + \alpha)(C'_1 + \alpha)}\bigg|$$

$$\leq b^2 \alpha^2 \bigg(\frac{1}{a} - \frac{1}{b}\bigg).$$

Hence

$$\frac{1}{2}\mathbf{E}\bigg[\bigg(\frac{C_1 \alpha}{C_1 + \alpha} - \frac{C'_1 \alpha}{C'_1 + \alpha}\bigg)^2\bigg] \leq b^4 \bigg(\frac{1}{b} - \frac{1}{a}\bigg)^2 \frac{1}{2}\mathbf{E}[(C_{n,1} + C_{n,2})^4].$$

Since both $C_{n,1}$ and $C_{n,2}$ are distributed as $C_{n-1}/2$, by Lemma 1, we have $C_{n,1} + C_{n,2} \leq 1/(a(n-1))$, and this yields

(4)
$$\frac{1}{2}\mathbf{E}\left[\left(\frac{C_1\alpha}{C_1+\alpha} - \frac{C_1'\alpha}{C_1'+\alpha}\right)^2\right] \leq \frac{1}{2}\left(\frac{b}{a}\right)^4\left(\frac{1}{b}-\frac{1}{a}\right)^2\frac{1}{(n-1)^4}$$
$$\stackrel{\text{def}}{=} \frac{K_0}{(n-1)^4}.$$

We now use the first term on the right-hand side of (3) to devise a recurrence. We have

$$\left|\frac{C_1(C_{n,1}+C_{n,2})}{C_1+C_{n,1}+C_{n,2}} - \frac{C_1(C_{n,1}'+C_{n,2})}{C_1+C_{n,1}'+C_{n,2}}\right|$$
$$= \frac{C_1^2|C_{n,1}-C_{n,1}'|}{(C_1+C_{n,1}+C_{n,2})(C_1+C_{n,1}'+C_{n,2})}$$
$$\leq |C_{n,1}-C_{n,1}'|.$$

Accordingly,

$$\mathbf{E}\left[\left(\frac{C_1(C_{n,1}+C_{n,2})}{C_1+C_{n,1}+C_{n,2}} - \frac{C_1(C_{n,1}'+C_{n,2})}{C_1+C_{n,1}'+C_{n,2}}\right)^2\right] \leq \mathbf{E}[(C_{n,1}-C_{n,1}')^2]$$
$$= 2\cdot\mathbf{Var}[C_{n,1}]$$
$$= \frac{1}{2}\cdot\mathbf{Var}[C_{n-1}].$$

Therefore, letting $K_1 = \max\{K_0, \mathbf{Var}[C_1]\}$ and recalling (4) we have the following recurrence relation:

$$\mathbf{Var}[C_n] \leq \frac{K_1}{(n-1)^4} + \frac{1}{2}\cdot\mathbf{Var}[C_{n-1}].$$

Expanding the recurrence yields readily

$$\mathbf{Var}[C_n] \leq K_1\cdot\sum_{i=1}^{n-1}\frac{2^{1-i}}{(n-i)^4} + \frac{1}{2^{n-1}}\mathbf{Var}[C_1] \leq K\cdot\sum_{i=1}^{n-1}\frac{2^{1-i}}{(n-i)^4} + \frac{K}{2^{n-1}}.$$

Since

$$\sum_{i=1}^{n-1}\frac{2^{1-i}}{(n-i)^4} + \frac{1}{2^{n-1}} \leq \sum_{i=1}^{\lfloor n/2 \rfloor}\frac{2^{1-i}}{(n-i)^4} + \sum_{i=\lfloor n/2 \rfloor+1}^{n-1}\frac{2^{1-i}}{(n-i)^4} + \frac{1}{2^{n-1}}$$
$$\leq \frac{2^5}{n^4}\sum_{i\geq 1}2^{-i} + 2^{1-\lfloor n/2 \rfloor}\sum_{i\geq 1}2^{-i} \leq \frac{2^{10}}{n^4},$$

the claimed bound on the variance holds as long as $K \geq K_1$. Finally, for all $n$, by Lemma 1 we have

$$
\mathbf{Var}[R_n] \leq \mathbf{E}\left[\left(R_n - \frac{1}{\mathbf{E}[C_n]}\right)^2\right] = \mathbf{E}\left[\left(\frac{\mathbf{E}[C_n] - C_n}{C_n \mathbf{E}[C_n]}\right)^2\right]
$$
(5)
$$
\leq b^4 n^4 \mathbf{Var}[C_n],
$$

so letting $K = \max\{b^4, 1\} \cdot K_1$, the proof is complete. □

We remark that this theorem is tight up to a constant factor unless $X$ is deterministic (in which case $\mathbf{Var}[C_n] = 0$). Indeed, by considering equation (2), since $C_{n,1}$ and $C_{n,2}$ are both of order $n^{-1}$, we see that fluctuations of constant size in the value of $C_1$ change $C_n$ by order $n^{-2}$. Such fluctuations occur with positive probability, so we must have that $\mathbf{Var}[C_n] \geq \varepsilon n^{-4}$ for some $\varepsilon > 0$ depending on $X$.

By Chebyshev's inequality, (5) also implies tail bounds on the probability that the resistance $R_n$ deviates from $1/\mathbf{E}C_n$.

COROLLARY 4. *There exists a constant $K$, such that for all $t > 0$ and $n \geq 1$, we have*

$$
\left\{\left|R_n - \frac{1}{\mathbf{E}C_n}\right| > t\right\} \leq \frac{K}{t^2}.
$$

*It follows that $\mathbf{E}|R_n - 1/\mathbf{E}C_n| \leq 1 + K$.*

PROOF. The first claim is immediate using (5) and Chebyshev's inequality; to see the second claim, observe that

$$
\mathbf{E}|R_n - 1/\mathbf{E}C_n| \leq 1 + \int_1^\infty \{|R_n - 1/\mathbf{E}C_n| \geq x\}\, dx \leq 1 + K \int_1^\infty \frac{dx}{x^2} = 1 + K.
$$
□

**3. The expected resistance and conductance.** In this section we give precise locations for the expected values $\mathbf{E}C_n$ and $\mathbf{E}R_n$, respectively. Let $\sigma^2 = \mathbf{Var}[X]$ and let $\mu = \mathbf{E}X$.

THEOREM 5. *There exist constants $M_1$ and $M_2$ depending only on $a$ and $b$ such that for all integers $n \geq 2$,*

$$
\left|\mathbf{E}R_n - \mu n + \frac{\sigma^2}{\mu} \ln n\right| \leq M_1 \quad \text{and} \quad \left|\mathbf{E}C_n - \frac{1}{\mu n} - \frac{\sigma^2 \ln n}{\mu^3 n^2}\right| \leq \frac{M_2}{n^2}.
$$

We remark that since $\mathbf{Var}[R_n]$ is certainly bounded from below by a positive constant (unless $X$ is deterministic, in which case we know $R_n$ precisely), we have determined the value of $\mathbf{E}R_n$ up to the order of its standard deviation. Furthermore, since $\mathbf{Var}[C_n]$ is of order $n^{-4}$, we have likewise determined $\mathbf{E}C_n$ up to the order of its standard deviation.

The techniques we use to handle the recurrence relation have been used by de Bruijn [8] to analyze slowly converging sequences and by Flajolet and Odlyzko [12] in the context of heights of simple trees.

PROOF OF THEOREM 5. We focus on $\mathbf{E}C_n$. By Corollary 4, bounds on $\mathbf{E}C_n$ immediately yield bounds on $\mathbf{E}R_n$. As in the proof of Theorem 2, we decompose $T_{n+1}$ into three independent conductors $C_1$, $C_{n+1,1}$ and $C_{n+1,2}$ (see Figure 2). Let $C_n$ and $C'_n$ be independent copies of the conductance between the root and level $n$. Since
$$C_{n+1} = \frac{C_1 \cdot (C_{n+1,1} + C_{n+1,2})}{C_1 + C_{n+1,1} + C_{n+1,2}},$$
$C_{n+1,1}$ and $C_{n+1,2}$ are both distributed as $C_n/2$, and $C_1$ is distributed as $1/X$, we have, in distribution,

$$(6) \qquad C_{n+1} = \frac{C_n + C'_n}{2} \cdot \frac{1}{1 + X((C_n + C'_n)/2)},$$

where $X$ is independent of all the other random variables appearing in (6). The second factor in (6) can be rewritten as

$$(7) \quad \frac{1}{1 + X((C_n + C'_n)/2)} = 1 - X\left(\frac{C_n + C'_n}{2}\right) + X^2 \cdot \left(\frac{C_n + C'_n}{2}\right)^2 \\ - X^3 \cdot \left(\frac{C_n + C'_n}{2}\right)^3 \cdot \frac{1}{1 + X((C_n + C'_n)/2)}.$$

Using the equality (7) to replace the term $1/(1 + X(C_n + C'_n)/2)$ in (6) and taking expectations, we obtain

$$(8) \quad \begin{aligned} \mathbf{E}C_{n+1} &= \mathbf{E}C_n - \frac{\mathbf{E}X}{2} \cdot (\mathbf{E}[C_n^2] + [\mathbf{E}C_n]^2) \\ &\quad + \frac{\mathbf{E}[X^2]}{4} \cdot (\mathbf{E}[C_n^3] + 3\mathbf{E}[C_n^2]\mathbf{E}C_n) \\ &\quad - \mathbf{E}\left[\frac{X^3(C_n + C'_n)^4}{16(1 + X((C_n + C'_n)/2))}\right], \end{aligned}$$

where we have used the equalities $\mathbf{E}[(C_n + C'_n)^2] = 2(\mathbf{E}[C_n^2] + [\mathbf{E}C_n]^2)$ and $\mathbf{E}[(C_n + C'_n)^3] = 2(\mathbf{E}[C_n^3] + 3\mathbf{E}[C_n^2]\mathbf{E}C_n)$. By Lemma 1, we have deterministically
$$\frac{a^3}{b^4 n^4} \cdot \frac{1}{1 + b/(an)} \leq \frac{X^3(C_n + C'_n)^4}{16(1 + X((C_n + C'_n)/2))} \leq \frac{b^3}{a^4 n^4},$$

so (8) yields

$$\mathbf{E}[C_{n+1}] = \mathbf{E}[C_n] - \frac{\mathbf{E}X}{2} \cdot (\mathbf{E}[C_n^2] + [\mathbf{E}C_n]^2)$$
(9)
$$+ \frac{\mathbf{E}[X^2]}{4} \cdot (\mathbf{E}[C_n^3] + 3\mathbf{E}[C_n^2]\mathbf{E}C_n) + O(n^{-4}),$$

where the order notation $O(\cdot)$ depends only on $a$ and $b$. We observe that, by Theorem 2,

(10) $\quad \mathbf{E}[C_n^2] + [\mathbf{E}C_n]^2 = \mathbf{Var}[C_n] + 2[\mathbf{E}C_n]^2 = 2[\mathbf{E}C_n]^2 + O(n^{-4}).$

Furthermore, since $\mathbf{E}[(C_n - \mathbf{E}C_n)^3] = O(n^{-1}) \cdot \mathbf{Var}[C_n] = O(n^{-5})$ by Theorem 2, we have

$$\begin{aligned}\mathbf{E}[C_n^3] &= \mathbf{E}[(C_n - \mathbf{E}C_n)^3] + 3\mathbf{E}[C_n^2]\mathbf{E}C_n - 3[\mathbf{E}C_n]^3 + [\mathbf{E}C_n]^3 \\ &= O(n^{-5}) + 3(\mathbf{Var}[C_n] + [\mathbf{E}C_n]^2)\mathbf{E}C_n - 2[\mathbf{E}C_n]^3 \\ &= 3\mathbf{Var}[C_n]\mathbf{E}C_n + [\mathbf{E}C_n]^3 + O(n^{-5}) \\ &= \mathbf{E}[C_n]^3 + O(n^{-5}),\end{aligned}$$

so

$$\begin{aligned}\mathbf{E}[C_n^3] + 3\mathbf{E}[C_n^2]\mathbf{E}C_n &= 4\mathbf{E}[C_n]^3 + 6\mathbf{Var}[C_n]\mathbf{E}C_n + O(n^{-5}) \\ &= 4\mathbf{E}[C_n]^3 + O(n^{-5}).\end{aligned}$$

Combining (9), (10) and (11), we obtain

$$\mathbf{E}C_{n+1} = \mathbf{E}C_n - \mathbf{E}X[\mathbf{E}C_n]^2 + \mathbf{E}X^2[\mathbf{E}C_n]^3 + O(n^{-4}).$$

Dividing through by $\mathbf{E}C_{n+1}\mathbf{E}C_n$ and letting $x_n = 1/\mathbf{E}C_n$ gives

(11) $\quad x_n = x_{n+1} - \mathbf{E}X \cdot \dfrac{x_{n+1}}{x_n} + \mathbf{E}[X^2] \cdot \dfrac{x_{n+1}}{x_n^2} + O(n^{-2}).$

We let $\delta_n = x_{n+1}/x_n - 1$ and let $\varepsilon_n = x_{n+1}/x_n^2$, and remark that $\delta_n$ and $\varepsilon_n$ are both $O(n^{-1})$. Summing (11) gives

(12) $\quad x_{n+1} = n\mathbf{E}X + \mathbf{E}X \cdot \sum_{i=1}^{n} \delta_i - \mathbf{E}[X^2] \sum_{i=1}^{n} \varepsilon_i + O(1).$

Since both $\delta_i$ and $\varepsilon_i$ are $O(i^{-1})$, (12) immediately yields the bound

(13) $\quad x_{n+1} = n\mathbf{E}X + O(\ln n) = (n+1)\mathbf{E}X + O(\ln n),$

a bound we will bootstrap to prove the theorem. From (11) we have

$$\frac{x_{n+1}}{x_n} = 1 + \frac{\mathbf{E}X + \mathbf{E}X \cdot \delta_n - \mathbf{E}[X^2] \cdot \varepsilon_n}{x_n} + O(n^{-3}),$$

so since $x_{n+1}/x_n$ also equals $1 + \delta_n$, solving for $\delta_n$ we obtain

$$\delta_n = \frac{\mathbf{E}X - \mathbf{E}[X^2]\varepsilon_n}{x_n - \mathbf{E}X} + O(n^{-3}) = \frac{\mathbf{E}X}{x_n - \mathbf{E}X} + O(n^{-2}) \tag{14}$$

as long as $n$ is large enough to ensure that $x_n - \mathbf{E}X$ does not happen to be zero (say $n \geq n_0$ for some fixed $n_0$ depending only on $a$ and $b$). Similarly,

$$\varepsilon_n = \frac{1}{x_n} \cdot \frac{x_{n+1}}{x_n} = \frac{1}{x_n} + \frac{\delta_n}{x_n} = \frac{1}{x_n} + O(n^{-2}) \tag{15}$$

for all $n \geq 1$. Combining (12), (14) and (15) gives the identity

$$\begin{aligned}
x_{n+1} - n\mathbf{E}X &= [\mathbf{E}X]^2 \sum_{i=n_0}^{n} \frac{1}{x_i - \mathbf{E}X} - \mathbf{E}[X^2] \sum_{i=n_0}^{n} \frac{1}{x_i} + O(1) \\
&= [\mathbf{E}X]^3 \sum_{i=n_0}^{n} \frac{1}{x_i(x_i - \mathbf{E}X)} - \mathbf{Var}[X] \sum_{i=n_0}^{n} \frac{1}{x_i} + O(1) \\
&= -\mathbf{Var}[X] \sum_{i=n_0}^{n} \frac{1}{x_i} + O(1).
\end{aligned} \tag{16}$$

Since $x_i = i\mathbf{E}X + O(\ln i)$ by (13), we have

$$\sum_{i=n_0}^{n} \frac{1}{x_i} = \sum_{i=n_0}^{n} \frac{1}{i\mathbf{E}X + O(\ln i)} = \sum_{i=n_0}^{n} \left( \frac{1}{i\mathbf{E}X} + \frac{O(\ln i)}{(i\mathbf{E}X)^2} \right) = \frac{\ln n}{\mathbf{E}X} + O(1),$$

so (16) yields

$$x_{n+1} - n\mathbf{E}X = \frac{\mathbf{Var}[X]}{\mathbf{E}X} \ln n + O(1).$$

The first assertion of the theorem follows immediately, and second assertion of the theorem follows since $x_{n+1} = \mathbf{E}R_{n+1} + O(1)$ by Corollary 4. □

**4. Sub-Gaussian tails bounds for the resistance.** In this section we show that the resistance $R_n$ does not only have a bounded variance but all its moments are also bounded and satisfies a sub-Gaussian tail inequality. In order to show this, we use a strengthening of the Efron–Stein inequality developed by Boucheron et al. [7], together with Thomson's formula. This flow-based formulation of the resistance was used by Benjamini and Rossignol [5] to show submean variance bounds for the random resistance in $\mathbb{Z}^2$. Given a graph $G$, let $E(G)$ be the set of edges of $G$. Recall that if $F$ denotes the set of unit flows from the root $r$ to depth $n$ in $T_n$, then

$$R_n = \inf_{\Theta \in F} \left\{ \sum_{e \in E(T_n)} r_e \Theta(e)^2 \right\}. \tag{17}$$

Furthermore, there is a unique unit flow $\Theta^*$ which attains the above infimum. As observed by Benjamini and Rossignol [5], it is a straightforward consequence of the Efron–Stein inequality that

$$\textbf{Var}[R_n] \leq \frac{(b-a)^2}{2} \sum_{e \in E(T_n)} \mathbf{E}[\Theta^*(e)^4]. \tag{18}$$

We now describe the result we use from [7] and how it can be combined with Thomson's formula to obtain a sub-Gaussian tail bound for $R_n$.

Suppose we are given independent random variables $\mathbf{U} = (U_1, \ldots, U_m)$ and a real-valued function $Z = f(U_1, \ldots, U_m)$. For integers $i = 1, \ldots, m$, let $U_i'$ be an independent copy of $u$, and let $Z_i' = f(U_1, \ldots, U_{i-1}, U_i', U_{i+1}, \ldots, U_m)$. Let

$$V^+ = \sum_{i=1}^{m} (Z - Z_i')_+^2 \quad \text{and let} \quad V^- = \sum_{i=1}^{m} (Z - Z_i')_-^2,$$

where $(\cdot)_+$ and $(\cdot)_-$ denote the positive and negative parts, respectively. Boucheron et al. [7] prove that if there exists a constant $C$ such that $V_+ \leq C$ almost surely, then

$$\{Z > \mathbf{E}Z + t\} \leq e^{-t^2/4C},$$

and if $V_- \leq C$ almost surely, then

$$\{Z < \mathbf{E}Z - t\} \leq e^{-t^2/4C}.$$

We shall apply these bounds with $\mathbf{U} = (X_e)_{e \in E(T_n)}$ and $Z = R_n$. From now on, $X_e'$ denotes an independent copy of $X_e$, and $R_n^{(e)'}$ denotes the resistance of $T_n$ when $X_e$ is replaced by $X_e'$ while all other resistances are kept unchanged. If $\Theta^*$ denotes the unit flow attaining the infimum in the expression of $R_n$ by Thomson's formula, and $\Theta^{*,e}$ is the minimizing unit flow when $X_e$ is replaced by $X_e'$, then by Thomson's formula and the deterministic bound $|X_e - X_e'| \leq 2^{d(e)}(b-a)$,

$$(R_n - R_n^{(e)})_+ \leq (X_e - X_e')_+ \cdot \Theta^{*,e}(e)^2,$$
$$\leq (b-a)2^{d(e)}\Theta^{*,e}(e)^2$$

and similarly,

$$(R_n - R_n^{(e)})_- \leq (b-a)2^{d(e)}\Theta^*(e)^2.$$

Thus,

$$V_+ \leq (b-a)^2 \sum_{e \in E(T_n)} 2^{2d(e)}\Theta^{*,e}(e)^4 \quad \text{and} \quad V_- \leq (b-a)^2 \sum_{e \in E(T_n)} 2^{2d(e)}\Theta^*(e)^4.$$

The key argument is the following deterministic bound.

LEMMA 6. *The optimal unit flow of any edge $e \in E(T_n)$ satisfies, deterministically,*

$$\Theta^*(e) \leq \frac{bn}{a(n-d(e)+1)2^{d(e)-1}}.$$

PROOF. Let $v$ denote the endpoint of $e \in E(T_n)$ closer to the root $r$ of $T_n$. If $U(v)$ is the voltage at vertex $v$ when the unit current $\Theta^*$ flows from the root $r$ to the leaves and the leaves have voltage 0, then by Ohm's law,

$$R_{n,e}\Theta^*(e) = U(v) \leq U(r) = R_n \leq bn,$$

where $R_{n,e}$ is the effective resistance of the subtree rooted at $v$. On the other hand, by Rayleigh's monotonicity law (see [9], page 53) and since all $X_e$'s are at least $a$, we have

$$R_{n,e} \geq a(n-d(e)+1)2^{d(e)-1}.$$

Comparing the two bounds, we obtain the bound of the lemma. □

Since the upper bound of Lemma 6 does not depend on the random values of the resistances $X_e$, the same inequality holds for $\Theta^{*,e}(e)$ as well. Thus, we have

$$V_+ \leq \frac{2^4 b^4 (b-a)^2}{a^4} \sum_{e \in E(T_n)} \left(\frac{n}{n-d(e)+1}\right)^4 2^{-2d(e)}$$

$$= \frac{2^4 b^4 (b-a)^2}{a^4} \sum_{i=1}^{n} \left(\frac{n}{n+1-i}\right)^4 2^{-i}$$

and the same upper bound applies to $V_-$ as well. Since the sum on the right-hand side is bounded, there exists a constant $C = C(a,b)$ such that both $V_+ \leq C$ and $V_- \leq C$. As a consequence, we have the following sub-Gaussian concentration inequality.

THEOREM 7. *There exists a constant $C$ depending on $a$ and $b$ only such that for every $t > 0$,*

$$\{|R_n - \mathbf{E}R_n| > t\} \leq 2e^{-t^2/4C}.$$

**5. Concluding remarks.** We conclude by discussing some possible extensions to branching random networks. When discussing the possible analogues of our results for branching processes, the following formulation of branching processes is useful. We start from a single "root edge" $uv$ and let $v$ be the root of a supercritical branching process with branching distribution $B$ that satisfies $\{B = 0\} = 0$. We use $\mathcal{T}$ to refer to this edge-rooted branching process. We say that a node $w \neq u$ has depth $i$ if there are $i$ edges on the path from $v$ to $w$ ($v$ has depth 0). For $i \geq 0$, we let $Z_i$ be the number of nodes of $\mathcal{T}$ at depth $i$—so in particular $Z_0 = 1$.

A branching random network is simply an edge-rooted branching process $\mathcal{T}$ as above. To each edge $e$ at depth $d$, we assign a random resistance $r_e = [\mathbf{E}B]^d X_e$, where the $X_e$ are independent, identically distributed positive random variables taking values in $[a,b]$ as in the binary case. As before, we let $C_n$ (resp., $R_n$) be the effective conductance (resp., effective resistance) from the root to depth $n$. As in the binary case, the above scaling most naturally corresponds to the critical case of a random walk in a random environment on branching processes with push-back. In particular, if $B$ is deterministically 2 then we recover the model of the previous sections.

It is easily seen that $C_n$ is not concentrated. Let $B_1$ be the number of children of the root—then as in the case of binary branching, we may first decompose $T_n$ into independent conductors $C_1$ and $C_{n,1}, \ldots, C_{n,B_1}$, so that

$$C_n = \frac{C_1(C_{n,1} + \cdots + C_{n,B_1})}{C_1 + C_{n,1} + \cdots + C_{n,B_1}}.$$

If $C_n$ is concentrated, then $C_{n,1} + \cdots + C_{n,B_1}$ is well approximated by $B_1 \cdot \mathbf{E}C_{n-1}/\mathbf{E}B_1$, so $C_n$ is close to

$$\frac{C_1 \cdot B_1 \cdot \mathbf{E}C_{n-1}/\mathbf{E}B_1}{C_1 + B_1 \cdot \mathbf{E}C_{n-1}/\mathbf{E}B_1} = \frac{\mathbf{E}C_{n-1}}{\mathbf{E}B_1} \cdot \left(\frac{1}{B_1} + \frac{\mathbf{E}C_{n-1}}{C_1 \cdot \mathbf{E}B_1}\right)^{-1}.$$

But the latter expression is *not* concentrated—a constant change in $B_1$ changes this expression by a constant factor. This should not be surprising: if the root has many offspring, the conductance is likely to be much (a constant factor) higher than if the root has a single child. It seems likely that at least the first-order behavior of $C_n$ and $R_n$ is governed by $W = \lim_{n \to \infty} Z_n/[\mathbf{E}B]^n$. If the resistance random variable $X$ is constant (say $X = 1$) then this is easily seen: the series-parallel laws give $R_n = \sum_{i=0}^n [\mathbf{E}B]^i/Z_i$, and $\mathbf{E}B^i/Z_i$ tends to $1/W$ $\mathcal{T}$-a.s., so $R_n/n$ tends to $1/W$ $\mathcal{T}$-a.s. If we additionally assume that $B$ has finite variance then this convergence is also in expectation [2], Theorem I.6.2, and it immediately follows that $C_n/\mathbf{E}C_n$ tends to $W/\mathbf{E}W$ a.s. and in expectation. In particular, this implies that $\lim_{n \to \infty} C_n/\mathbf{E}C_n$ has absolutely continuous distribution (as long as $B$ is not constant; see [2], Theorem I.10.4), which is the "scaled analogue" of Question 4.1 from Lyons, Pemantle and Peres [16] mentioned in the Introduction.

We would expect that even when $X$ is not constant, for *any* $\lambda$ with $1 \leq \lambda \leq \mathbf{E}B$, in the network where the resistances of edges at depth $i$ are scaled by $\lambda^i$, $\lim_{n \to \infty} C_n/\mathbf{E}C_n$ has an absolutely continuous distribution. In the special case that $\lambda = \mathbf{E}B$, we would venture that $R_n/n$ tends to $\mathbf{E}X/W$, $\mathcal{T}$-a.s. and in expectation. However, we were unable to extend the arguments

EFFECTIVE RESISTANCE OF RANDOM TREES 15

used to prove Theorems 2 and 5 to the branching process. In particular, once we condition on $\mathcal{T}$, our techniques for manipulating (3) and (6) in order to devise recurrences fail, most notably because in this setting the identical distribution of subtrees at equal depth is lost. It seems plausible that as in the case of binary branching $R_n - W \cdot (n\mathbf{E}X)$ is $O(\ln n)$, again $\mathcal{T}$-a.s. and in expectation, but the coefficient of $\ln n$ also seems likely to depend on $\mathcal{T}$ and the precise nature of this dependence is unclear to us.

Finally, observe that Barndorff-Nielsen [3] and Barndorff-Nielsen and Koudou [4] have noticed an interesting link between inverse Gaussian (or reciprocal inverse Gaussian) random variables and effective resistances of random networks: if resistances of the edges are distributed like i.i.d. inverse Gaussian random variables, then the effective resistance of the entire tree is distributed like a reciprocal inverse Gaussian.

**Acknowledgments.** We wish to thank an anonymous referee for many comments and corrections, and in particular for substantially simplifying the proof of Theorem 7.

L. Addario-Berry
Département de mathematiques
  et de statistique
Université de Montréal
C.P. 6138, succ. centre-ville
Montréal, QC, H3C 3J7
Canada
E-mail: addario@dms.umontreal.ca

N. Broutin
Projet Algorithmes
INRIA Rocquencourt
78153 Le Chesnay
France
E-mail: nicolas.broutin@inria.fr

G. Lugosi
ICREA and Department of Economics
Pompeu Fabra University
Ramon Trias Fargas 25-27
08005, Barcelona
Spain
E-mail: gabor.lugosi@gmail.com